\documentclass[10pt]{article}

\usepackage{amsfonts}
\usepackage[utf8x]{inputenc}
\usepackage[T1]{fontenc}
\usepackage{lmodern}
\usepackage{amsmath,amsthm,amssymb,bm}
\usepackage[hidelinks]{hyperref}
\usepackage{cleveref,nameref}
\usepackage[a4paper]{geometry}
\usepackage{mathtools}
\usepackage{enumitem}

\usepackage[doublespacing]{setspace}

\newtheorem{theo}{Theorem}
\newtheorem{lem}{Lemma}[section]

\newtheorem{coro}{Corollary}[theo]
\newtheorem{defi}{Definition}
\newtheorem{propri}{Property}
\newtheorem{rem}{Remark}

\DeclareMathOperator{\sign}{sign}

\newcommand{\modint}{\displaystyle\copy\tratto\kern-10.4pt\int\limits}

\def\Mu{{M^1}}
\def\u{\mu}

\def\R{\mathbb R}
\def\N{\mathbb N}

\DeclareMathOperator{\supp}{supp}
\DeclareMathOperator*{\essinf}{ess\,inf}

\def\OV{\overline}

\def\O{\Omega}

\def\Max{\mathop{\rm Max\,}}
\def\infess{\mathop{\rm inf\,ess}}
\def\supess{\mathop{\rm sup\,ess}}
\def\a{\alpha}
\def\b{\beta}
\def\d{\delta}
\def\eps{\varepsilon}
\def\f{\varphi}

\def\p{\partial}
\def\om{\omega}
\def\s{\sigma}

\def\LEQ{\leqslant}
\def\GEQ{\geqslant}

\def\tee{{\theta\eta}}
\def\si{\hbox{{if}\ }}

\def\DST{\displaystyle}
\def\pr{{\bf Proof :\\ }}
\def\HF{\hfill{$\diamondsuit$\\ }}

\def\div{{\rm div\,}}
\def\meas{{\rm measure\,}}
\def\exp{{\rm exp\,}}

 \def\NN{\nonumber}

\def\Log{{\rm log\,}}

\def\dist{{\rm dist\,}}

\def\calD{\mathcal D}

\def\OV{\overline}
\def\ORA{\overrightarrow}
\def\OVO{{\overline \Omega}}
\def\BSL{\backslash}

\def\PRI{{^{\prime\prime}}}

\def\Vinf{{V,\infty}}

\def\Wkq{{W^{k,q}}}
\def\BkK{{{\mathbb B}^k_K}}
\def\BdK{{{\mathbb B}^2_K}}

\def\xj{\xrightarrow[j\to+\infty]{}0}
\def\xji{\xrightarrow[j\to+\infty]{}}
\def\Cap{{\rm Cap}}
\def\calO{{\mathcal O}}

\begin{document}

\title{ \Large \bf  Potential-capacity and some applications}

\author{
	{\bf  Jean Michel Rakotoson}\\%
		{\small Universit\'e de Poitiers,}\\
	{\small Laboratoire de Math\'ematiques et Applications }\\
	{\small  UMR CNRS 7348 - SP2MI, France.}\\	
	{\small Bat H3 - Bd Marie et Pierre Curie,  T\'el\'eport 2,}\\
	{\small  F-86962 Chasseneuil Futuroscope Cedex, France.}\\
	{\small e-mail : rako@math.univ-poitiers.fr}}
\date{}
\maketitle

\begin{abstract}
We introduce a  new capacity associated to a non negative function~$V$. We apply this notion  to the study of a necessary and sufficient condition to ensure the existence and uniqueness of a Schr\"odinger   type equation with  measure data and with  an operator whose coefficients are discontinuous. Namely, for a  potential $V$, $f$ a bounded Radon measure on 
$\O$, then the equation $\cal L$$_Vu=-\Delta u+U\cdot\nabla u+Vu=f$ has a solution in $\DST L^1(V) \cap L^1_0(\O)=\Big\{g\hbox{ measurable }, \int_\O|g|Vdx\hbox{ is finite and }
\lim_{\eps\to0}\int_{\!\!\{x:\dist(x;\partial\O)\LEQ\eps\}}\!\!\!\!\!\!\!\!\!\!\!\!\!\!\!\!\!\!\!\!\!\!\!\!\!\!\!\!\!|g|dx=0\ \Big\}$
if and only if f does not charge "irregular points" of $V$, provided that the set of "irregular points" have a zero potential capacity. As a byproduct of our results, we have the non existence of a Green operator for some $\cal L$$_V$.\\
Our method is also based on a new topology  and density of $C^2_c(\O\BSL K)$ in $C^2_0(\OV\O)$ whenever $K$ has a zero potential-capacity.
\end{abstract}
\textbf{Keywords :}  {\it Capacity, Schr\"odinger equations, potential, measures.}

\textbf{AMS Classification} 35J75, 35J15, 35J25, 35J10, 76M23

\section{Introduction}
In recent works (see \cite{DGRT}, \cite{DGR}), we have studied the equation 
$$-\Delta \omega + U \cdot \nabla \omega+ V\omega = f$$ in  a smooth open bounded domain $\Omega$ whenever the potential $V$ is locally integrable on the domain, under the  Dirichlet condition $u=0$ on $\partial\O$.\\
The first natural question is: what happens if we remove this integrability condition  on $V$?\\
When examining the prototype of $V$ say $V(x)= |x-a|^{-m}$, with $m>0$, $ a \in \O,$
we observe that there is an interaction between the point $a$, the power $m$ and the right hand side $f$.\\ To describe the growth of $V$ and such interaction, we introduce here a new capacity associated to the potential $V$. Roughly speaking, the more $V$ contains "irregular points" the more its capacity will be  small. In particular, we will focus on potential whose "irregular points" are of capacity  zero.
\\This new capacity  is slightly different to the usual one considered by many authors (see \cite{APonce} \cite{Ziemer}) for a complete review). Indeed, we recall that, if $K$ is a compact subset of an open set $\O$ of $\R^n$, then,\\ for $1\LEQ k<+\infty,\ 1\LEQ q<+\infty$ the $\Wkq$ capacity of $K$ is usually defined as
\begin{equation}\label{eq1}
\Cap_\Wkq(K)=\inf\Big\{||\f||^q_{\Wkq(\O)},\ \f\in \BkK\big\}
\end{equation}
where
\begin{equation}\label{eq2}
\BkK=\Big\{\psi\in C^k_c(\O),\ 0\LEQ\psi\LEQ1,\ \psi=1 \hbox{ in a neighborhood of }K\Big\}.
\end{equation}

Here, we shall consider a potential $V\GEQ0$ on $\O,\ V\not \equiv0$ for $\psi\in C^2_c(\O)$, we define
\begin{equation}\label{eq3}
||\psi||_{V,\infty}=||\psi||_{L^1(\O)}+\left\|\dfrac{\nabla\psi}{\sqrt V}\right\|_\infty+\left\|\dfrac{\Delta\psi}{V}\right\|_\infty,
\end{equation}
and we shall associate, the following capacity function, for a compact $K$ included in $\O$
\begin{equation}\label{eq4}
\Cap_{V,\infty}(K)=\inf\Big\{||\psi||_{V,\infty},\ \psi\in{\mathbb B}^2_K\Big\}
\end{equation}

Such capacity possesses common properties as for the above classical capacities (see Section \ref{s1} below), namely, we will show in particular that $$\hbox{if  $\Cap_\Vinf(K_i)=0,\ i\in J$(finite) then }\DST\Cap_\Vinf\Big(\bigcup_{i\in J}K_i\Big)=0.$$ 
Roughly speaking, such capacity will measure how singular is the potential $V$? And how "large" is this singularity. For instance, if $a\in\O$ and $V$ behaves like $ |x-a|^{-m}$ near $a$, then $\Cap_\Vinf(\{a\})=0$ if $m\GEQ2$ and $\Cap_\Vinf(\{a\})>0\ \si m<2$. 
But one of the most important properties that we need for the applications  are  :

\begin{theo}\label{t1}\ \\
Let $K$ be compact included in $\O$. Assume that $\Cap_{V,\infty}(K)=0$. Then
there exists a sequence $(\psi_j)_j$,\ \  $\psi_j\in C^2_c(\O)$ such that
\begin{enumerate}

\item $\psi_j(x)\DST\xrightarrow[j\to+\infty]{}0$ for a.e in $\O$ and strongly in $L^1(\O)$, $\psi_j=1$ on $K$.\\
$$\hbox{$\dfrac{|\nabla\psi_j|}{\sqrt V}\DST\xrightarrow[j\to+\infty]{}0$ and
$\dfrac{\Delta\psi_j}V\DST\xrightarrow[j\to+\infty]{}0$ strongly in $L^\infty(\O).$}$$

\item  If furthermore 
$ V\in L^{\frac n2,1}_{loc}(\O_K)$  with $\O_K=\Big\{x\in\O,\ \dist(x;K)>0\Big\}$, then,\\ for all $x\in \O-K$
$$\psi_j(x)\xrightarrow[j\to+\infty]{}0,$$
more precisely, if $\O_{K,0}\subset\subset\O_K$, then
$$\Max_{\OVO_{K,0}}|\psi_j(x)|\xrightarrow[j\to+\infty]{}0,\quad||\nabla\psi_j||_{L^{n,1}(\O_{K,0})}\xj.$$

\item If $a\in\O,\ V(x)=|x-a|^{-m}$ then $V$ is in $L^{\frac n2,1}(\O)$ if and only if $m<2.$.
\end{enumerate}
\end{theo}

The natural question is then, can we give   sufficient conditions to ensure that
\begin{equation}\label{eq5}
\Cap_{V,\infty}(K)=0?
\end{equation}
The answer to that question is naturally linked with the motivations of our study.
 One of them  is the following :\\
 
 { Let $\mu_0$ be the Dirac mass at the origin, $m$ a positive parameter, then we  observe the following phenomena : \\
 If $m\GEQ2$ then there is no solution of $$(M_1)\qquad\begin{cases}-\Delta u(x)+\dfrac{u(x)}{|x|^m}=\mu_0&\hbox{in }B(0;1)\subset\R^n,\ n\GEQ2,\\
 u(x)=0&\hbox{if }|x|=1.\end{cases}$$
 But if $m<2$, the above problem $(M_1)$ possesses at least one solution $u$. The same phenomena were also given in \cite{BenBr}. \\
 Let us notice that $(M_1)$ has a solution if $n
=1$.}\\

Another motivation that we shall prove in this note is the following removable type singularities result :

{\bf Proposition} {\bf (removable singularities with potential)} \label{p1}\ \\
{\it 
Assume (for simplicity) that $V(x)=\DST\sum_{i=1}^p\dfrac{b_i}{|x-a_i|^{m_i}},\ m_i\GEQ0$, $b_i>0,\ a_i\in\O$ and consider $K=\Big\{a_i,\ m_i\GEQ2\}$,
$w\in L^1(\O;V)\cap L^1_0(\O))$ such that $\forall\, \f\in C^2_c(\O\backslash K)$ we have
$$\int_\O w(-\Delta\f+V\f)dx=0.$$ 
Then
$$w\equiv0.$$
Here $\d(x)={\rm distance}(x;\p\O)$.
} 
\\

So, the natural question is that if we consider an arbitrary potential $V\GEQ0$, how can we replace the set of singularities $K=\Big\{a_i,\ m_i\GEQ2\Big\}$? 

The question seems to be linked with some density problem ( with an adequate topology).\\
The tough problem linked with that question is the construction of an appropriate sequence smooth function vanishing over $K$ and disappearing when we pass to the limit for an adequate topology. These are the purpose of our main results stated in the next section. Namely a generalization of the above proposition for a large class of potential $V$ and applications to some existence and non existence result for weak or very weak solution. We shall provide few examples of compact $K$ whose $(V,\infty)$-capacity is zero.

\section{Notations Definitions - Primary definitions and  results}\label{s1}
We shall keep the notation we used to employ. 
We set 
$$
	L^0(\O)=\Big\{v:\O\to\R\hbox{ Lebesgue measurable}\Big\}
$$ 
and we denote    by $L^p(\O)$ the usual Lebesgue space $1\LEQ p\LEQ +\infty$.
Although it is not too often used, we shall use the notation 
$$W^{1,p}(\O)=W^1L^p(\O)$$
for the  associated Sobolev space. We need the following definitions :
\begin{defi}{\bf of the distribution function and monotone rearrangement}\ \\
Let $u\in L^0(\O)$. The distribution function of $u$ is the decreasing function
	\begin{eqnarray*}
		m=m_u:\R &\to &[0,|\O|] \\
			t & \mapsto & \meas\big\{ x:u(x)>t\big\}=|\big\{u>t\big\}|.
	\end{eqnarray*} 
	The generalized inverse $u_*$ of $m$ is defined by, for $s\in[0,|\O|[$,
	$$
		u_*(s)=\inf\Big\{t:|\big\{u>t\big\}|
	\LEQ s\Big\},
	$$
	and is called the decreasing rearrangement of $u$. We shall set $\O_*=]0,|\O|\,[.$
\end{defi}
\begin{defi}\label{d5}\ \\
	Let $1\LEQ p\LEQ +\infty,\ 0<q\LEQ+\infty$ :
	\begin{itemize}
		\item If $q<+\infty$, one defines the following norm for $u\in L^0(\O)$
	$$\|u\|_{p,q}=\|u\|_{L^{p,q}}:=\left[\dfrac1{|\O|}\int_{\O_*}\left[t^{\frac1p}|u|_{**}(t)\right]^q\frac{dt}t\right]^{\frac1q}\hbox{ where }|u|_{**}(t)=\dfrac1t\int_0^t|u|_*(\s)d\s.$$
		\item If $q=+\infty$,
	$$\|u\|_{p,\infty}=\sup_{0<t\LEQ|\O|}t^{\frac1p}|u|_{**}(t).$$
	The space 
	\begin{equation} 
		L^{p,q}(\O)=\Big\{ u\in L^0(\O):\|u\|_{p,q}<+\infty\Big\}
	\end{equation}
	 is called a {\bf Lorentz space}.
	
		\item If $p=q=+\infty,$ $L^{\infty,\infty}(\O)=L^\infty(\O).$	
		\item The dual of $L^{1,1}(\O)$ is called $L_\exp(\O)$
	\end{itemize}
\end{defi}

\begin{rem}\ \\
We recall that $L^{p,q}(\O)\subset L^{p,p}(\O)=L^p(\O)$ for any $p>1,\ q\GEQ1$.
\end{rem}

\begin{defi}\ \\
If $X$  is a Banach space in $L^0(\O)$, we shall denote the Sobolev space associated to $X$ by
$$W^1X=\Big\{\f\in L^1(\O):\nabla\f\in X^n\Big\}$$ or more generally for $m\GEQ1$,
$$ W^mX=\Big\{\f\in W^1X,\ \forall\,\a=(\a_1,\ldots,\a_n)\in\N^n,\ |\a|=\a_1+\ldots+\a_n\LEQ m,\ D^{|\a|}\f\in X\Big\}.$$
We also set
$$W^1_0X=W^1X\cap W_0^{1,1}(\O).$$
\end{defi}

We also need to recall the Hardy's inequality in $L^{n^{\prime },\infty }$ \
saying that if $\O$ is a bounded Lipschitz domain :\\

\begin{equation}
\int_{\Omega }\left|\frac{u }{\delta }\right|^q\LEQ c\Vert \nabla u \Vert^q
_{L^{n^{\prime },\infty }}\qquad \forall u \in W_{0}^{1}L^{n^{\prime
},\infty }(\Omega ),
\end{equation}
with $n' = \frac n {n-1}$, $1\LEQ q<n'$. This inequality can be obtained from the results of \cite{RakoJFA} (see also \cite{DR}) since $W_0^1 L^{n',\infty} (\Omega) \subset W_0^1 (\Omega; 1 + |\log \delta|)$.\\

We need the following Lemma whose proof is given in \cite{Ziemer, Kufner, Triebel}
\begin{lem}\label{lA1}\ \\
Let $A\subset\R^n$ be closed and for $x\in\R^n$ let $d(x)=d(x;A)$ denote the distance from $x$ to $A$. Let 
$$U=\Big\{x:d(x)<1\Big\}.$$
Then there is a function $\rho\in C^\infty(U-A)$ and a positive number $M=M(n)$ such that
$$M^{-1}d(x)\LEQ\rho(x)\LEQ M\,d(x),\ x\in U-A$$
$$|D^\a\rho(x)|\LEQ c(\a)\,d(x)^{1-|\a]},\ x\in U-A,\ \  |\a|=\a_1+\ldots+\a_n.$$
In particular, the result holds if $A=\partial\O$ boundary of an open bounded set $\O$ , in this case $$\rho\in C^\infty(\O)\hbox{ and }d(x)=\d(x)=\dist(x;\p\O).$$
\end{lem}
\begin{defi}{\bf of $(\Vinf)$-capacity or potential-capacity}\\
Let $V\GEQ0$ be a measurable function on $\O$, $V$ non identically zero, V is called a potential function.\\
The $(\Vinf)$-capacity of a compact $K$ included in $\O$  (or potential-capacity of $K$) is given by relation (\ref{eq4}).
\end{defi}
We will denote by $c$ different constant, sometimes we will specify the dependence with respect to the data.

\begin{propri}{\bf of $(\Vinf)$-capacity}\\
For any compact $K$ in $\O$, we have
\begin{enumerate}

\item measure$(K)\LEQ \Cap_\Vinf(K)$.

\item If $K_1 $ is another compact included in $K$ then
$$\Cap_\Vinf(K_1)\LEQ \Cap_\Vinf(K).$$

\item For all $\eps>0$, there exists an open set $\om$ containing $K$ such that for all compact $K'$ satisfying $K\subset K'\subset\om$, on has
$$\Cap_\Vinf(K')\LEQ \Cap_\Vinf(K)+\eps.$$

\item If $V_1,\ V_2$ are two nonnegative potential $V_1\LEQ V_2$ then
$$\Cap_{V_2,\infty}(K)\LEQ \Cap_{V_1,\infty}(K).$$
\end{enumerate}
\end{propri}

\pr
\begin{enumerate}

\item For $\psi\in\BdK$ we have  $\DST\meas(K)\LEQ\int_\O\psi(x)dx\LEQ||\psi||_\Vinf$ which gives the result.

\item If $K_1\subset K $ then $\BdK\subset{\BdK}_{\!_1}$. Therefore
$$\Cap_\Vinf(K_1)\subset \Cap_\Vinf(K).$$

\item Let $\eps>0$, then there exists $\psi_\eps\in\BdK$ such that
\begin{equation}\label{eq6}
||\psi_\eps||_\Vinf -\eps\LEQ \Cap_\Vinf(K)\LEQ||\psi_\eps||_\Vinf.
\end{equation}
Since $\psi_\eps=1$ in a neighborhood of $K$, thus there exists an open set of $\om$ on which $\psi_\eps=1$. Then for all compact $K'$ with 
$K\subset K'\subset\om$ one has $\psi_\eps=1$ on $K' $ and then $\psi_\eps\in\BdK_{'}\ $. Thus,
$$\Cap_\Vinf(K')\LEQ||\psi_\eps||_\Vinf\LEQ \Cap_\Vinf(K)+\eps.$$

\item If $0\LEQ V_1\LEQ V_2$ then $\dfrac1{V^\a_2}\LEQ\dfrac1{V^\a_1}$ if $\a=\dfrac12, \a=1$ from which we get the result. ${\  }$\qquad \HF

\end{enumerate}

\begin{rem}\ \\
\begin{itemize}
\item In the definition of $(\Vinf)$-capacity, we can add a  different power on the potential $V$ but the choice of the power is linked with the applications.

\item The property (3) is the so-called continuity from the right in Choquet's capacity theory.
\end{itemize}
\end{rem}

{\bf Proof of Theorem \ref{t1}}
\begin{enumerate}

\item As for relation (\ref{eq6}) considering $\eps=\dfrac1j,\ j\GEQ1$ we have a sequence $(\psi_j)_j$:
\begin{equation}\label{eq7}
\psi_j=1\hbox{ on } K,\qquad 0\LEQ||\psi_j||_\Vinf-\Cap_\Vinf(K)\LEQ\dfrac1j.
\end{equation}
If $\Cap_\Vinf(K)=0$ then $$||\psi_j||_\Vinf\xrightarrow[j\to+\infty]{}0$$
which implies
$$\left\|\dfrac{\nabla\psi_j}{\sqrt V}\right\|_\infty\xj,\ \left\|\dfrac{\Delta\psi_j}V\right\|_\infty\xj,\hbox{ and }||\psi_j||_{L^1}\xj.$$
This last convergence implies that for a subsequence still denoted by $\psi_j$ that $$\psi_j(x)\xj\hbox{ for a.e.}$$

\item Let $x\in \O_K$, there exists $r>0$ so that $B(x;r)\subset \O_K$. From Poincar\'e-Sobolev's inequality or P.D.E. regularity (see \cite{FRR})
\begin{eqnarray*}
||\nabla \psi_j||_{L^{n,1}(B(x;r))}&\LEQ&
c\left[\,||\psi_j||_{L^{\frac n2,1}(B(x;r)\,)}+||\Delta\psi_j||_{L^{\frac n2,1}(B(x;r)\,)}\right]\\
&\LEQ&c\left[\,||\psi_j||_{L^1(\O)}+\left\|\dfrac{\Delta\psi_j}V\right\|_\infty||V||_{L^{\frac n2,1}(B(x;r)\,)}\right]\xj.
\end{eqnarray*}

\begin{eqnarray*}
\Max_{y\in\OV B(x;r)}\left|\psi_j(y)\right|
&\LEQ& c||\nabla\psi_j||_{L^{n,1}(B(x;r))}+c||\psi_j||_{L^1(\O)}\xj.
\end{eqnarray*}
If $\O_{K,0}$ is open set relatively compact in $\O$ by recovering $\OVO_{K,0}$ and using the same argument as the above result we deduce
$$\Max_{y\in\OVO_{K,0}}|\psi_j(y)|\xj,\qquad||\nabla\psi_j||_{L^{n,1}(\O_{K,0})}\xj.$$
\item A direct and simple computation shows that
$$||\,|x-a|^m||_{L^{\frac n2,1}}<+\infty\hbox{ if and only if } m<2.$$
\ $\qquad $\HF
\end{enumerate}

\section{Few examples of compact $K$ having a $(\Vinf)$-capacity zero}
\begin{theo}{\bf (Comparison near a compact)}\label{t2}\\
Let $K$be a compact in $\O$, $V_1$ and $V_2$ two nonnegative potentials satisfying
\begin{enumerate}

\item $\exists\eta>0$ such that $V_1\LEQ V_2$ on a compact set $$K_\eta=\Big\{x\in\O,\ d(x;K)\dot=\dist(x;K)\LEQ2\eta\Big\}\subset\O.$$

\item $V_1 $ is bounded from below and above on  
$$\calO_\eta=\Big\{x\in\O:\dfrac12\eta<d(x;K)<2\eta\Big\}\hbox{ i.e }0<\infess_{\calO_\eta }V_1\LEQ\supess_{\calO_\eta}V_1<+\infty.$$
\end{enumerate}
Then there exists a constant $c_\eta>0$ such that 
\begin{equation}\label{eq8}
\Cap_{V_2,\infty}(K)\LEQ c_\eta\Cap_{V_1,\infty}(K).
\end{equation}
In particular, 
\begin{equation}\label{eq9}
\hbox{if }\Cap_{V_1,\infty}(K)=0\hbox{ then }\Cap_{V_2,\infty}(K)=0.
\end{equation}
\end{theo}
\pr
Let $\theta$ be in $C^\infty_c(\O)$ such that $0\LEQ\theta\LEQ1$ and 
$\theta=1$ on $\Big\{x\in\O:d(x,K)\LEQ~\eta~\Big\}$\\ $\hbox{ and support}(\theta)\dot=\supp(\theta)\subset\Big\{x\in\O:d(x;K)<\dfrac32~\eta\Big\}$.\\
Let us show that there exists a constant $c_\theta>0$ such that  for all $\psi\in \BdK,$
\begin{equation}\label{eq10}
||\theta\psi||_{V_2,\infty}\LEQ c_\theta||\psi||_{V_1,\infty}
\end{equation}
We need the following
\begin{lem}\label{l1}\ \\
There exists a constant $c_{\theta\eta}>0$ such that for all $\psi\in \BdK$ we have
\begin{equation}\label{eq11}
|\psi(x)|\LEQ c_{\theta\eta}\Big[||\psi||_{L^1(\O)}+||\nabla\psi||_{L^\infty(\calO_\eta)}\Big] \quad\forall\,x\in\supp(\theta)\cap\Big\{\eta\LEQ d(\cdot;K)\LEQ\dfrac32\eta\Big\}.
\end{equation}
\end{lem}

\pr
By the compactness of the set 
$$H_0=\supp(\theta)\cap\Big\{\eta\LEQ d(\cdot;K)\LEQ\dfrac32\eta\Big\}$$
we have a family $(B(x_i;r_i))_{i=1,\ldots,p} $ such that 
$$H_0\subset\bigcup_{i=1}^pB(x_i;r_i)=\calO\subset\calO_\eta.$$
Applying the Sobolev embedding, we have a constant $c_{\theta\eta}>0$
$$||\psi||_{L^\infty(\calO)}\LEQ c_{\theta\eta}\Big[||\psi||_{L^1(\O)}+||\nabla\psi||_{L^\infty(\calO_\eta)}\Big],\quad \psi\in\BdK.$$
This gives the result\HF

Let $\psi\in\BdK$ then $||\theta\psi||_{L^1}\LEQ||\psi||_{L^1}$. For $x\in\O$ 
$$\nabla(\theta\psi)(x)=\nabla\theta\psi(x)+\theta(x)\nabla\psi(x).$$
We distinguish 3 cases
\begin{enumerate}

\item If $x\notin\supp(\theta) $ then $\nabla(\theta\psi)(x)=0$, therefore we have
$$\dfrac{\nabla(\theta\psi)}{\sqrt V_2}(x)=0\LEQ||\psi||_{V_1,\infty}.$$

\item If $x\in\supp\theta,\
,\ x\notin\calO_\eta$ then $V_1(x)\LEQ V_2(x)$ and 
$$|\nabla(\theta\psi)(x)|\LEQ|\nabla\psi(x)|:\hbox{ so that }\left|\dfrac{\nabla(\theta\psi)}{\sqrt V_2}(x)\right|\LEQ\left|\dfrac{\nabla\psi}{\sqrt V_1}(x)\right|\LEQ||\psi||_{V_1,\infty}.$$

\item If $x\in\supp\theta,\ x\in\calO_\eta$ we still have $V_1(x)\LEQ V_2(x)$ but $\DST V_1(x)\GEQ\essinf_{\calO_\eta} V_1 >0$ so that
\begin{enumerate}

\item if $d(x;K)\LEQ\eta,\ \nabla(\theta\psi)(x)=\nabla\psi(x)$ so we still have
$$\left|\dfrac{\nabla(\theta\psi) }{\sqrt V_2}(x)\right|\LEQ\left|\dfrac{\nabla\psi}{\sqrt V_1}(x)\right|\LEQ||\psi||_{V_1,\infty},$$

\item if $\eta<d(x;K)\LEQ\dfrac32\eta$, we use Lemma \ref{l1} to 
$$\left|\dfrac{\nabla(\theta\psi)}{\sqrt V_2}(x)\right|\LEQ c_\theta\left[|\psi(x)|+\dfrac{|\nabla\psi|}{\sqrt V_1}(x)\right]\LEQ c_{\theta\eta}||\psi||_{V_{1,\infty}}.$$
\end{enumerate}
\end{enumerate}

Since $\Delta(\theta\psi)(x)=\Delta\theta(x)\psi(x)+2\nabla\theta(x)\nabla\psi(x)+\theta(x)\Delta\psi(x)$, we can argue as before to deduce that

\begin{equation}\label{eq12}
\left|\dfrac{\Delta(\theta\psi)}{\sqrt V_2}(x)\right|\LEQ c_\tee||\psi||_{V_{1,\infty}},\ \forall\,x\in\O.
\end{equation}

We have shown
\begin{equation}\label{eq13}
\left\|\dfrac{\nabla(\theta\psi)}{\sqrt V_2}\right\|_{L^\infty(\O)}+\left\|\dfrac{\Delta(\theta\psi)}{V_2}\right\|_{L^\infty(\O)}\LEQ c_\tee||\psi||_{V_{1,\infty}}.
\end{equation}
Thus we deduce relation (\ref{eq10})  $\forall\,\psi\in\BdK$.  We then have
$$\Cap_{V_{2,\infty}}(K)\LEQ c_\theta\Cap_{V_{1,\infty}}(K).$$
\ \qquad\HF

Here are few examples of Theorem \ref{t2}.
\begin{coro}{\bf of Theorem \ref{t2}}\\
Let $A$, be a closed set included in $\O$ whose measure is zero, $m>2,\ m\in\R,\ V$ a  potential such that\\ there exists $\eta>0,\ c>0$ with $V(x)\GEQ\dfrac c{d(x;A)^m}$ for $x\in \Big\{y:d(y;A)\LEQ2\eta\Big\}$. Then
$$\Cap_\Vinf(A)=0.$$
\end{coro}
\pr
Let us set $V_1(x)=\dfrac c{d(x;A)^m},\ x\in\O$. According to Theorem \ref{t2} it is sufficient to show that$$\Cap_{V_{1,\infty}}(A)=0.$$
Let $H\in C^\infty(\R)$ such that
\begin{equation}\label{eq2000}H\in C^\infty(\R)\hbox{ such }H(t)=\begin{cases}1&\si t\GEQ2,\\0&\si t\LEQ1.\end{cases}\end{equation}
  and denote by $\d(x)=\dist(x;\p\O)$.\\
According to Lemma \ref{lA1} that we have a function $\rho\in C^\infty(\O)$, two constants $c_1>0,\ c_2>0$ such that 
\begin{enumerate}

\item  $c_1\d(x)\LEQ\rho\LEQ c_2\d(x),\ \forall\,x\in\O$

\item $\forall\a=(\a_1,\ldots,\a_n)\in \N^n,\ \exists c_\a>0$ such that $$|D^\a\rho(x)|\LEQ c_\a\rho(x)^{1-|\a|}\hbox{ with }
|\a|=\a_1+\ldots+\a_n,\ 
D^\a=\dfrac{D^{\a_1+\ldots+\a_n}}{\p x_1^{\a_1}\ldots\p x_n^{\a_n}}.$$
 More, we have $ \rho_A\in C^\infty(\O\BSL A),\ M=M(n)>0,\ c(\a),\ |\a|\LEQ2$ such that for all $x\in\O\BSL A,$
 
 \item $M^{-1}d(x;A)\LEQ \rho_A(x)\LEQ M\,d(x;A)$,
 
 \item $|D^\a \rho_A(x)|\LEQ c(\a)d(x;A)^{1-|\a|}.$
\end{enumerate}
Consider the sequence $\psi_j(x)=\big(1-H(j \rho_A(x))\big)H\big(j\rho(x)\big)$. Then $\psi_j\in C^\infty_c(\O)$ and $j\GEQ j_a$, large enough so that $\Big\{x\in\O: \rho_A(x)<\dfrac1{j_a}\Big\}\subset\Big\{ x\in\O:\dist(x;A)<\dist(A;\partial\O)\Big\}$.
 $$\psi_j(x)=\begin{cases}
 1&\si  \rho_A(x)<\dfrac1j,\\
 0&\si  \rho_A(x)\,\GEQ\, \dfrac2j,\\
 1-H\big(j \rho_A(x)\big)&\si \dfrac1j< \rho_A(x)<\dfrac2j.\end{cases}$$
 On the set $D_j=\Big\{\dfrac1j< \rho_A(x)<\dfrac2j\Big\}$ one has
 $$\nabla\psi_j(x)=-jH'(j \rho_A(x))\nabla{ \rho_A(x)},$$ so that
\begin{equation}\label{eq15}
 \rho_A(x)\,|\nabla\psi_j(x)|\LEQ c_1||H'||_\infty j \rho_A(x)\LEQ c_{1H}
\end{equation}

and
$$\Delta\psi_j(x)=-H^\PRI(j \rho_A(x))j^2|\nabla\, \rho_A(x)|^2-H'(j \rho_A(x))j\Delta{ \rho_A(x)}.$$

From which we have
\begin{equation}\label{eq16}
 \rho_A(x)^2|\Delta\psi_j(x)|\LEQ c_3||H\PRI||_\infty(j \rho_A(x))^2+||H'||_\infty c_4j \rho_A(x)\LEQ c_{2H}.
\end{equation}
Since the measure of $A$ is zero and $\psi_j(x)\xj\quad\forall\,x\in\O\BSL A,$\\
we deduce by the Lebesgue dominated theorem that
$$||\psi_j||_{L^1}\xj.$$
Since $\Delta\psi_j(x)=\nabla\psi_j(x)=0$ outside of $D_j$, we deduce from the above estimates  
\begin{equation}\label{eq17}
\left\|\dfrac{\nabla\psi_j}{\sqrt V_1}\right\|_{L^\infty(\O)}\LEQ c'_{1H}\dfrac1{j^{\frac m2-1}}\xj
\end{equation}
and
\begin{equation}\label{eq18}
\left\|\dfrac{\Delta\psi_j}{V_1}\right\|_{L^\infty(\O)}\LEQ c'_{2H}\dfrac 1{j^{m-2}}\xj.
\end{equation}
Since $\psi_j\in{\mathbb B}^2_A$, we deduce
$$\Cap_{V_1,\infty}(A)\LEQ||\psi_j||_{V_1,\infty}\xj.$$
\ \qquad\HF 


\begin{coro}{\bf of Theorem \ref{t2}}\\
Let $S_1=\Big\{x\in\R^n:|x|=1\Big\}$ the unit sphere of $\R^n$, $m>2$ assume that $S_1\subset\O$ and let $V$ a nonnegative potential such that there exists $\eta>0,\ c>0$ such that $V(x)\GEQ\dfrac c{|\,|x|-1|^m}$ for all $x\in\{y\in\O,\ d(y;S_1)\LEQ2\eta\Big\}$. Then
$$\Cap_\Vinf(S_1)=0.$$
\end{coro}
One important property concerns the potential-capacity of a finite union of compact $\DST\bigcup_{i\in J}K_i$ such that $\Cap_\Vinf(K_i)=0$ we are not able to prove the subadditivity, but we also have:
\begin{theo}\label{t3}\ \\
Let $V$ be a nonnegative potential, $K_i,\ i\in J$ be a finite number of compact sets included in $\O$. Assume that $\Cap_\Vinf(K_i)=0\ \forall\,i\in J$. Then
$$\Cap_\Vinf(\bigcup_{i\in J}K_i)=0.$$
\end{theo}
\pr
Since $\Cap_\Vinf(K_i)=0$ there exists a sequence $\psi_{ij}\in C^2_c(\O)$ such that for a.e $x$, $$\psi_{ij}(x)\xj,\quad||\psi_{ij}||_\Vinf\xj,$$ $\psi_{ij}(x)=1$ in a neighborhood of $K_i,\ 0\LEQ\psi_{ij}\LEQ1$. Let us consider $H\in C^\infty(\R), \  0\LEQ H\LEQ1$ as in relation (\ref {eq2000}), $ \rho\in C^\infty(\O)$ equivalent to the distance function $\d(x)=\dist(x;\p\O)$.\\
Since $\psi_{ij}\in C^2_c(\O)$ then we have a set $\Delta_{ij}\subset\Big\{x\in\OVO:\psi_{ij}(x)=0\Big\}$ which is an open neighborhood of the boundary. Therefore, we can consider the open set $\Delta_j=\bigcap_{i\in J}\Delta_{ij}$ neighborhood of $\p\O$. \\
Since $\a_j=\dist(\O\backslash\Delta_j;\p\O)>0$, we can consider a sequence $\mu_j>0$, such that $\mu_j<\a_j$ and $\mu_j\to0$ as $j\to+\infty$. One has, in this case, the set $$\Big\{x\in\O:\rho(x)\LEQ\mu_j\Big\}\subset\Delta_j,$$ 
otherwise, we will have a point $x$ such $\rho(x)\LEQ\mu_j$ and $x\in\O\backslash\Delta_j$ so  that $$\dist(\O\backslash\Delta_j;\p\O)\LEQ\rho(x).$$

The function $$\Phi_j(x)=\big(1-\prod_{i\in J}(1-\psi_{ij}(x))\big)H\big(\dfrac2{\mu_j}\rho(x)\big)\hbox{ with }3\mu_j<\dist\big(\bigcup_{i\in J}K_i;\p\O\big)$$ satisfies
 \begin{enumerate}
\item $\DST\Phi_j(x)=1,\ x\in \bigcup_{i\in J} K_i$,
\item $\Phi_j\in C^2_c(\O)$,
\item $0\LEQ \Phi_j(x)\LEQ1,\ \ \Phi_j(x)=1-\DST\prod_{i\in J}\big(1-\psi_{ij}(x)\big)\ \si\rho(x)>\mu_j,\ \ ~{\Phi_j(x)=0\ \si\rho(x)\LEQ~\mu_j.}$
\end{enumerate}

We shall set for simplicity $\Phi_{ij}(x)=1-\psi_{ij}(x)$.\\
For $x\in\O$ such that $\rho(x)>\mu_j$, we have $H\Big(\dfrac2{\mu_j}\rho\Big)=1 $ and
\begin{equation}\label{eq19}
\nabla\Phi_j(x)=-\sum_{k\in J}\prod_{i\in J, k\neq i}\Phi_{ij}(x)\nabla\Phi_{kj}(x)
\end{equation}
 
\begin{equation}\label{eq20}
\left|\dfrac{\nabla\Phi_j}{\sqrt V}(x)\right|
\LEQ\sum_{k\in J}\left|\dfrac{\nabla \Phi_{kj}}{\sqrt V}(x)\right|
=\sum_{k\in J}\left|\dfrac{\nabla\psi_{kj}}{\sqrt V}(x)\right|
\LEQ\sum_{k\in J}||\psi_{kj}||_\Vinf.
\end{equation}
We also has
\begin{eqnarray}
|\Delta\Phi_j(x)|&\LEQ&\sum_{k\in J}|\Delta\psi_{kj}(x)|+\sum_{k\in J}\sum_{\ell\in J}|\nabla \psi_{kj}(x)|\,|\nabla\psi_{\ell j}(x)|\NN\\
\left|\dfrac{\Delta\Phi_j}V(x)\right|
&\LEQ&\sum_{k\in J}\left|\dfrac{\Delta\psi_{kj}}V(x)\right|
+\sum_{(k,\ell)\in J^2}\left|\dfrac{\nabla\psi_{kj}}{\sqrt V}(x)\right|\,
\left|\dfrac{\nabla \psi  _{\ell j}}{\sqrt V}(x)\right|\NN\\
&\LEQ&\sum_{k\in J}||\psi_{kj}||_\Vinf
+\sum_{(k,\ell)\in J^2}||\psi_{kj}||_\Vinf||\psi_{\ell j}||_\Vinf.\label{eq21}
\end{eqnarray}

If $\rho(x)\LEQ\mu_j$ then $x\in \Delta_j$ and   
$$1-\prod_{i\in J}(1-\psi_{ij}(x))=0:\Phi_j(x)=0.$$
We conclude that relations (\ref{eq20}) and (\ref{eq21}) hold true. Therefore, we always have
\begin{equation}\label{eq22}
\left\|\dfrac{\nabla\Phi_j}{\sqrt V}\right\|_{L^\infty(\O)}\LEQ\sum_{k\in J}||\psi_{kj}||_\Vinf,
\end{equation}
\begin{equation}\label{eq23}
\left\|\dfrac{\Delta\Phi_j}V\right\|_{L^\infty(\O)}
\LEQ\sum_{k\in J}||\psi_{kj}||_\Vinf+\left(\sum_{k\in J}||\psi_{kj}||_\Vinf\right)^2.
\end{equation}
On other hand, by the Lebesgue dominated convergence theorem, we have
\begin{equation}\label{eq24}
||\Phi_j||_{L^1(\O)}\xj.
\end{equation}
Relations (\ref{eq22}) to (\ref{eq24}) yield that
$$||\Phi_j||_\Vinf\xj.$$
Since $$\Cap_\Vinf\Big(\bigcup_{i\in J} K_i\Big)\LEQ||\Phi_j||_\Vinf,$$
we derive the result.\HF

As a consequence of Theorem \ref{t3}, we have\\
\begin{coro}{\bf of Theorem \ref{t3}}\\
For $i\in\{1,\ldots,m\}$, let $a_i\in\O,\ r_i\GEQ0,\ m_i>2,\ c_i>0$ real numbers.\\ Define $S_i=\big\{x\in \R^n:|x-a_i|=r_i\big\},$ $K=\DST\bigcup_{i=1}^m S_i$ assumed to be included in $\O$. Let $V$ be a nonnegative potential such that there exists $\eta>0$ such that 
$$V(x)\GEQ\sum_{i=1}^m\dfrac{c_i}{(|x-a_i|-r_i)^{m_i}}\hbox{ on } 
\big\{y:\dist(y;K)\LEQ\eta\big\}.$$
Then
$$\Cap_\Vinf(K)=0.$$
\end{coro}
\pr
We have seen in Corollary 2.1  of Theorem \ref{t2} that $\Cap_\Vinf(S_i)=0$ whenever $S_i\subset\O$. Applying Theorem \ref{t3}, we deduce the result.\HF

In the above Corollary 1 and 2 of Theorem \ref{t2} we  may replace $S_1$ by any compact  included in $\O$ whose measure  is zero. Concrete examples for application are given in \cite{CKS, PT}.

As we have announced in the introduction,, we have $\Cap_{|x|^{-m},\infty}(\{0\})>0$ if $m<2$. Here is the proof
\begin{theo}\label{t80}\ \\
Let $V$ be a nonnegative potential, $a\in\O$ be such that there exist $\eta>0$, $c>0$
$$V(x)\LEQ\dfrac c{|x-a|^m},\quad x\in B(a;2\eta)\hbox{ for some }m<2.$$
Then
$$\Cap_{\Vinf}(\{a\})>0.$$
\end{theo}
\pr
Let us set  $V_1(x)=\dfrac c{|x-a|^m},\ x\in\O\BSL\{a\}$.\\ Following Theorem \ref{t2}, $\Cap_\Vinf(\{a\})\GEQ c_\eta \Cap_{V_1,\infty}(\{a\})$, $ c_\eta>0$.\\ We have for $\f\in{\mathbb B}^2_{\{a\}}$,
\begin{eqnarray*}
||\nabla\f||_{L^{n,1}(B(a;\eta))}
&\LEQ&\Big[\,||\f||_{L^{\frac n2,1}(B(a;\eta))}+||\Delta\f||_{L^{\frac n2,1}(B(a,\eta))}\Big]
\\&\LEQ& c||\f||_{V_1,\infty}\Big[1+||V_1||_{L^{\frac\eta2,1}(B(a;\eta))}\Big]
\\&\LEQ& c_1||\f||_{V_1,\infty}<+\infty.
\end{eqnarray*}
Applying the Sobolev-Lorentz embedding
$$||\f||_{L^\infty(B(a;\eta))}\LEQ c_2\Big[||\f||_{L^{n,1}(B(a;\eta))}+||\nabla\f||_{L^{n,1}(B(a;\eta))}\Big]\LEQ c_3||\f||_{V_1,\infty}.$$
Since $\f(x)=1$ in a neighborhood of $a$, this last inequality implies  $1\LEQ c_3\Cap_{V_1,\infty}(\{a\})$, this implies the result.
\HF
\begin{rem}\ \\
\begin{enumerate}

\item As we state before, the choice of the power $\dfrac12$ and 1 in the definition is linked with the application, it is clear we can use other power as $\left\|\dfrac{\nabla\psi}{V^\a}\right\|$ and $\left\|\dfrac{\Delta\psi}{V^\b}\right\|$, $\a>0,\ \b>0$ (see \cite{RakoBook}).

\item In Corollary 2.1 of Theorem 2, we may take $m=2$, but the proof to show that \\ $\Cap_\Vinf(A)=0$ uses a different argument ( \cite{RakoBook} work in progress)
\end{enumerate}
\end{rem}
We define 
$$C^2_0(\OVO)=\Big\{\f\in C^2(\OVO),\ \f=0\hbox{ on }\p\O\Big\}.$$

\section{Approximation of functions in $C^2_0(\OVO)$}
We shall introduce the following sets :
$$L^1(\O;V)\dot=L^1(V)=\Big\{g:\O\to\R\hbox{ measurable such that }\int_\O|g(x)|V(x)dx<+\infty\Big\}$$
$$L^1_0(\O)=\Big\{g\in L^1(\O)\hbox{ such that  }\lim_{\eps\to0}\dfrac1\eps\int_{\{x:\d(x)\LEQ\eps\}}|g(x)|dx=0\Big\}$$

\begin{rem}\ \\
One has
$$L^1(\O;\d)=\Big\{g\in L^1(\O):\int_\O|g(x)|\dfrac{dx}{\d(x)}<+\infty\Big\}\hbox{ is {\bf strictly} included in $L^1_0(\O)$.}$$ Indeed, it was shown in \cite{RakoJFA} that  $$\si f\GEQ0\quad f\in L^1(\O;\d)\BSL L^1\big(\O;\d(1+|\Log\d|)\big)$$ { then the unique solution }$g\in L^1_+(\O) $ of
$$-\int_\O g\Delta\f dx=\int_\O \f f dx\quad\forall\,\f\in C^2_0(\OVO)$$
verifies
$$\int_\O\dfrac g\d(x)dx=+\infty.$$
But A. Ponce (\cite{APonce}, chap.20) shows that we have $g\in L^1_0(\O).$
\end{rem}
\ 

\begin{defi}\ \\
Let $\phi$ be in $C^2_0(\OVO)$. We will say that a sequence $(\f_j)_j$ of $C^2_0(\OVO)$ converges weakly in the  sense of the potential $V$ to $\phi$ if for all $g\in L^1(V)\cap L^1_0(\O):$
\def\xji{\DST\xrightarrow[j\to+\infty]{}}
\begin{enumerate}

\item $\DST\int_\O g\f_j Vdx\xji \int_\O g\,\phi\, Vdx$

\item $\DST\int_\O  g\p_k\f_j dx\xji \int_\O g\p_k\phi dx,\ \ k=1,\ldots,n$

\item $\DST\int_\O g\Delta\f_j\xji\int_\O g\Delta\phi dx.$
\end{enumerate}
Here $\p_k$ is the partial derivative with respect to the $k^{th}$ derivate.
\end{defi}
\ 
\begin{defi}\ \\
Let $\phi$ be  in $C^2_0(\OVO)$. We will say that a sequence $(\f_j)_j$ of $C^2_0(\OVO)$ converges weakly-strongly in the sense of $V$ to $\phi$ if  for all $g\in L^1(V)\cap L^1_0(\O)$
$$0=\lim_j\int_\O|g|\,|\f_j-\phi| Vdx=\lim_j\int_\O|g|\,|\nabla\f_j-\nabla\phi|dx=\lim_j\int_\O|g|\,|\Delta\f_j-\Delta\phi|dx.$$
\end{defi}
We have the
\begin{theo}\label{t4}\ \\
Let $K$ be a compact in $\O$ and $V$ a nonnegative potential. Assume that $\Cap_\Vinf(K)~=~0$. Then $$\hbox{the set $C^2_c(\O\BSL K)$ is \underline{\bf weakly-strongly} {\bf dense} in $C^2_0(\OVO)$ in the sense of the potential $V$.}$$
\end{theo}
\pr 
Let $\Phi$ be in $C^2_0(\OVO)$. Since $\Cap_\Vinf(K)=0$ one has a sequence $(\f_j)_j$,\ \ $\f_j\in C^2_c(\O)$ such that\\ $0\!\LEQ\!\f_j\!\LEQ\!\!1$,
$\f_j=1$ in a neighborhood of $K$ and  $||\f_j||_\Vinf\xj,\,\f_j\xj$ a.e in~$\O$. 
\\ As before, we then have a sequence $(\mu_j)_j$ tending to zero $\mu_j>0$ such that  $$\hbox{the set $\Big\{x:\d(x)\LEQ\mu_j\Big\}$ is included in $\Big\{x:\f_j(x)=0\}.$}$$
Let  $H$ be the function given in (\ref{eq2000}).\\
Then the sequence $\Phi_j=(1-\f_j)H\left(\dfrac2{\mu_j}\rho\right)\Phi$ where $\rho$ is the smooth function equivalent to the distance function $\d$ with the conditions $$||\nabla\rho||_\infty<+\infty,\ \  ||\rho\,\Delta\rho||_\infty<+\infty \hbox{ (see Lemma \ref{lA1})}$$ $\Phi_j$ belongs to $C^2_c(\O\BSL K)$. We have the following pointwise relations
\begin{equation}\label{eq30}
\Phi_j(x)\xji\Phi(x)\hbox{ a.e in }\O.
\end{equation}
If $\rho(x)>\mu_j,$  we have :
\begin{equation}\label{eq31}
\begin{matrix}
\ \ \Phi_j(x)=&\big(1-\f_j(x)\big)\Phi(x)\\
\nabla\Phi_j(x)=&-\nabla\f_j(x)\Phi(x)+\big(1-\f_j(x)\big)\nabla\Phi(x)\\
\Delta\Phi_j(x)=&-\Delta\f_j(x)\Phi(x)-2\nabla\f_j(x)\nabla\Phi(x)+\big(1-\f_j(x)\big)\Delta\Phi(x).\end{matrix}
\end{equation}
If $\rho(x)\LEQ\mu_j$, we know that $\f_j(x)=0$ so that
\begin{equation}\label{eq32}
\Phi_j(x)=H\left(\dfrac2{\mu_j}\rho(x)\right)\Phi(x)\ \si \ \rho(x)\GEQ\dfrac12{\mu_j}\hbox{ and }\Phi_j(x)=0\hbox{ otherwise.}
\end{equation}
Then on $\Big\{\dfrac12\mu_j\LEQ\rho\LEQ\mu_j\Big\}$, one has
$$\nabla\Phi_j(x)=\dfrac2{\mu_j}\nabla\rho H'\left(\dfrac2{\mu_j}\rho(x)\right)\Phi(x)+H\left(\dfrac2{\mu_j}\rho(x)\right)\nabla\Phi(x)$$
and

\begin{eqnarray}\label{eq33}
\Delta\Phi_j(x)&=&\dfrac2{\mu_j}\Delta\rho(x)\Phi(x)H'\left(\dfrac2{\mu_j}\rho(x)\right)+\left(\dfrac2{\mu_j}\right)^2\Phi(x)|\nabla\rho(x)|^2H^{\PRI}\left(\dfrac2{\mu_j}\rho(x)\right)\NN\\
&&+\dfrac4{\mu_j}H'\left(\dfrac2{\mu_j}\rho(x)\right)\nabla\rho(x)\cdot\nabla\Phi(x)+H\left(\frac2{\mu_j}\rho(x)\right)\Delta\Phi(x).
\end{eqnarray}
Let $g\in L^1(V)\cap L^1_0(\O)$. By the Lebesgue dominated theorem, we have
\begin{equation}\label{eq34}
\lim\int_\O|g(x)|\,|\Phi_j(x)-\Phi(x)|\,V(x)dx=0
\end{equation}
From relation (\ref{eq31}), we derive
\begin{eqnarray}\label{eq36}
\int_{\{\rho>\mu_j\}}|\nabla\Phi_j(x)-\nabla\Phi(x)|\,|g(x)|dx&\LEQ&
||\nabla\Phi||_\infty\int_\O|\f_j(x)|\,|g(x)|dx\\
&&+||\f_j||_\Vinf||\Phi||_\infty\left(\int_\O|g(x)|V(x)dx\right)\left(\int_\O|g|dx\right).\NN
\end{eqnarray}
(We have used the Cauchy Schwarz inequality : $\DST\int_\O|g|\sqrt Vdx\LEQ\left(\int_\O|g|Vdx\right)\left(\int_\O|g|dx\right)$\ )\\
Setting $A_j=\Big\{\dfrac12\mu_j\LEQ\rho\LEQ\mu_j\Big\}$, one has \begin{equation}\label{eq37}
\int_{\{\rho\LEQ\mu_j\}}|g(x)|\,|\nabla\Phi_j(x)-\nabla\Phi(x)|dx\LEQ
\int_{\{\rho\LEQ\mu_j\}}|g(x)|\,|\nabla\Phi(x)|dx+\int_{A_j}|g(x)|\,|\nabla\Phi_j(x)|dx.
\end{equation}
The first integral tends to zero using Lebesgue dominated theorem, while the second integral can be bound as
\begin{equation}\label{eq38}
\int_{A_j}|g(x)|\nabla\Phi_j(x)|dx\LEQ c_\Phi\left[\dfrac1{\mu_j}\int_{A_j}|g(x)|dx+\int_{A_j}|g(x)|dx\right]\xj\hbox{ since }g\in L^1_0(\O).
\end{equation}
From relations (\ref{eq36}) to (\ref{eq38}) we derive
\begin{equation}\label{eq39}
\lim\int_\O|g(x)|\,|\nabla\Phi_j(x)-\nabla\Phi(x)|dx=0.
\end{equation}
Using the same argument as above, we have
\begin{equation}\label{eq40}
\lim\int_\O|g(x)|\,|\Delta\Phi_j(x)-\Delta\Phi(x)|dx=0.
\end{equation}
Indeed, we have
\begin{eqnarray}\label{eq41}
\int_{\{\rho\GEQ\mu_j\}}|g(x)|\Delta \Phi_j(x)-\Delta\Phi(x)|dx
&\LEQ&\int_{\{\rho\GEQ\mu_j\}}|\Delta\f_j|\,|\Phi|\,|g|dx
+\int_{\{\rho\GEQ\mu_j\}}|\f_j|\,|\Delta\Phi|\,|g|dx\NN\\
&&+2\int_{\{\rho\GEQ\mu_j\}}|\nabla\f_j|\,|\nabla\Phi|\,|g|dx\NN\\
&\LEQ&c_\Phi||\f_j||_\Vinf\left[\int_\O|g|Vdx\right]\,\left[1+\int_\O|g|dx\right]\NN\\
&&+c_\Phi\int_\O|\f_j|\,|g|dx\xj.
\end{eqnarray}
On $\{\rho\LEQ\mu_j\},$ we have :
\begin{equation}\label{eq42}
\int_{\{\rho\LEQ\mu_j\}}|g|\,|\Delta\Phi_j-\Delta\Phi|dx\LEQ\int_{\{\rho\LEQ\mu_j\}}|g|\,|\Delta\Phi|+\int_{A_j}|g|\,|\Delta\Phi_j|dx.
\end{equation}
The first term tends to zero, while for the last term we replace $\Delta\Phi_j$ by its expression~:
\begin{equation}\label{eq43}
\int_{A_j}|g|\,|\Delta\Phi_j|dx\LEQ I_{1j}+I_{2j}+I_{3j}+I_{4j}.
\end{equation}
Using the fact that $\left|\dfrac{\Phi(x)}{\rho(x)}\right|\LEQ c||\nabla \Phi||_\infty$, and $|\rho\Delta\rho(x)|\LEQ c_2$ for all $x\in\O$, we have :
\begin{equation}\label{eq44}
I_{1j}\LEQ c\dfrac1{\mu_j}\int_{A_j}\left|\dfrac{\Phi(x)}{\rho(x)}\right|\,|\Delta\rho(x)\ \rho(x)|\,|g(x)|dx\LEQ c\dfrac1{\mu_j}\int_{A_j}|g(x)dx,
\end{equation}

\begin{equation}\label{eq46}
I_{2j}\LEQ c\left(\dfrac1{\mu_j}\right)^2\int_{A_j}|g(x)|\,\rho(x)dx\LEQ c\dfrac1{\mu_j}\int_{A_j}|g(x)|dx,
\end{equation}
\begin{equation}\label{eq46}
I_{3j}\LEQ c\dfrac1{\mu_j}\int_{A_j}|g(x)|dx,\qquad I_{4j}\LEQ c\int_{A_j}|g(x)|dx
\end{equation}
Thus
\begin{equation}\label{eq47}
\int_{A_j}|g|\,|\Delta\phi_j|dx\LEQ c\left[\dfrac1{\mu_j}\int_{A_j}|g(x)|dx+\int_{A_j}|g(x)|dx\right],
\end{equation}
the constant $c$ is independent of $j$ and $g$.
From relations (\ref{eq41}) to (\ref{eq47}), we derive the result.\HF

One may also give sufficient conditions to ensure that a sequence converges weakly in the sense of $V$.\\

Here is an example of such result :
\begin{theo}\label{t200}\ \\
Let $(\f_j)_j$ be a sequence of $C^2_0(\OVO)$, $K$ a compact in $\O,$ $V$ a nonnegative potential such that $V$ is upper semi-continuous, that is for all real $t$, the set $\{ V\GEQ t\}$ is closed in $\O$, and assume also that the set  $\Big\{x:V(x)=+\infty\Big\}$ is of measure zero, and:
\begin{enumerate}
\item $||\f_j||_{V,\infty}$ remains bounded in $\R_+$,

\item there exists $\f\in C^2_0(\OVO)$ such that the sequence $(\f_j)_j$ converges to $\f$ in $L^\infty(\O)$-weak-star.
\end{enumerate}
Then, $(\f_j)_j$ converges weakly to $\f$ in the sense of the potential $V$.
\end{theo}
{\bf Sketch of the proof}\\
Let $1>\eta>0$, the, $\Big\{x\in \O:V(x)\GEQ\dfrac1\eta\Big\}$ is closed in $\O$ thus $\O_\eta=\Big\{x\in\O: \dfrac1V(x)>\eta\Big\}$ is open and we have a constant $M$ such that $\forall\,j,\ \forall\,\eta\in]0,1[$
\begin{equation}\label{eq200}
||\nabla\f_j||_{L^\infty(\O_\eta)}+||\Delta\f_j||_{L^\infty(\O_\eta)}\LEQ\eta^{-1 } ||\f_j||_{V,\infty}\LEQ M\eta^{-1}.
\end{equation}
Then we deduce that for all $\eta\in]0,1[,\ \DST||\f_j-\f||_{C^1(\OVO_\eta)}\xrightarrow[j\to+\infty]{}0.$

Since  $\DST\O\BSL\bigcup_{\eta>0}\O_\eta$  is of measure zero, therefore,
$$\dfrac{\nabla \f_j}{\sqrt V}\rightharpoonup \dfrac{\nabla \f}{\sqrt V}\hbox{ and }\dfrac{\Delta\f_j}{\sqrt V}\rightharpoonup\dfrac{\Delta\f_j}V\hbox{ in $L^\infty$-weak-star when }j\to+\infty.$$
From those convergences, we derive the result.
\HF

{\bf An example of sequence satisfying Theorem \ref{t200}}\\
As example, we can take $A$ as in Corollary  2.1 of Theorem \ref{t2}, $V(x)=d(x;A)^{-2}$ and \\$\psi_j(x)=\Big(1-H\big(j\rho_A(x)\big)\Big)H\big(j\rho(x)\big)$, as in the proof of Corollary 2.1 of Theorem~\ref{t2}. Then, for $\f\in C^2_0(\OVO)$ the sequence $\f_j(x)=\big(1-\psi_j(x)\big)\f(x)$ satisfies conditions 1. and 2. .

Indeed, since $\f_j(x)\to \f(x)\ \forall\,x\in\O\BSL A$ and $||\f_j||_\infty\LEQ||\f||_\infty$, we deduce that $(\f_j)_j$ converges to $\f$ in $L^\infty(\O)$-weakly-star.\\
The set $\{x:V(x)=+\infty\}=A$ is of measure zero and $V$ is upper semi continuous.\\

More, $\f_j(x)=\f(x)$ if $\rho_A(x)\GEQ\dfrac2j$, and $\nabla\f_j(x)=0$ if $\rho_A(x)<\dfrac1j$,\\ on $D=\Big\{\dfrac1j<\rho_A<\dfrac2j\Big\}$, we have 
\begin{equation}\label{eq201}
\rho_A(x)|\nabla\f_j(x)|\LEQ c_{1\f}\hbox{ and }\rho_A(x)^2|\Delta\f_j(x)|\LEQ c_{2\f}.
\end{equation}
This implies $||\f_j||_{\Vinf}\LEQ M<+\infty.$\HF
\begin{coro}{\bf of Theorem \ref{t200}}\ \\
Let $V(x)=d(x;A)^{-2},\ A$ a compact set of measure zero in $\O$. Then
$$\hbox{ $C^2_c(\O\BSL A)$ is weakly dense in $C^2_0(\OVO)$ in the sense of the potential $V$.}$$
\end{coro}
\pr
Let $\f\in C^2_0(\OVO)$ then $\f_jH\big(j\rho(x)\big)$ is in $C^2_c(\O\BSL A)$ and  the above arguments imply the statement.\HF

\section{Applications of the potential-capacity and the approximation of $C^2_0(\OVO)$}
As a first application of the above results, we shall prove a removable type problem.\\

\begin{theo}\label{t5}\ \\
Let $K$ be compact included in $\O$. 
Assume that $C^2_c(\O\BSL K)=C^2_c(\O_K)$ is weakly dense in $C^2_0(\OVO)$ in the sense of potential $V$ and let $w\in L^1(\O;V)\cap L^1_0(\O)$ be such that  for all $\f\in C^2_c(\O_K)$ we have
\begin{equation}\label{eq51}
\int_\O w(-\Delta\f+V\f)dx=0.
\end{equation}
Then $w$ satisfies the same equation (\ref{eq51}) with $\f\in C^2_0(\OVO)$.
\end{theo}
\pr
Let $\f$ be in $C^2_0(\OVO)$. Then, we have a sequence $(\f_j)_j,\ \f_j\in C^2_c(\O\BSL K)$ such that

\begin{equation}\label{eq52}
\lim_{j\to_\infty}\int_\O w\f_j \, Vdx=\int_\O w\f\, V dx.
\end{equation}
\begin{equation}\label{eq53}
\lim_{j\to_\infty}\int_\O w\Delta\f_j dx=\int_\O w\Delta\f dx.
\end{equation}
Since $\DST0=\int_\O w(-\Delta\f_j+ V\f_j)dx$ thus we have the result by passing to the limit.~\HF

Next, we recall the following Kato's inequality (see \cite{APonce, MV, Kato, BP}).

\begin{lem}{\bf Kato's inequality and weak maximum principle}\label{l2}\\
Assume that $w$ and $f$ are in $L^1(\O)$ such $-\Delta w= f$ in $\calD'(\O)$. Then
\begin{enumerate}

\item \begin{equation}\label{eq54}
-\Delta|w|\LEQ f\sign(w)\hbox{ in }\calD'(\O),
\end{equation}

\item  $$-\Delta w_+\LEQ f\sign_+(w)\hbox{ in }\calD'(\O),$$

\item$$\si -\Delta w\LEQ0\hbox{ in } C^2_0(\OV\O)'\hbox{ (dual space) then }w\LEQ0.$$
\end{enumerate}
$$\sign_+(\s)=\begin{cases}1&\si\s>0,\\0&\hbox{otherwise},\end{cases}\hbox{\ \ and\ \ }
\sign(\s)=\begin{cases}1&\si\s>0,\\0&\si\s=0,\\-1&\si\s<0.\end{cases}$$
\end{lem}

\begin{coro}{\bf of Theorem \ref{t5} and Lemma \ref{l2}}\\
Under the same assumption as for Theorem \ref{t5}, the function $w$ verifying relation (\ref{eq51}) satisfies $$w\equiv0.$$
\end{coro}

\pr
Since $\calD(\O)=C^\infty_c(\O)\subset C^2_0(\OVO)$, then following Theorem \ref{t5}, we have
$$-\Delta w=-Vw\in L^1(\O)\hbox{ in }\calD'(\O).$$
From Kato's inequality, one has :
$$-\Delta(|w|)\LEQ-V|w|\LEQ0.$$
Therefore, using the same arguments as for Theorem \ref{t5},  the inequality holds  in the dual space $C^2_0(\OV\O)'$, we conclude that  $|w|\LEQ0:w=0$.\HF

Let $V$ be a nonnegative potential and define the  subset of $\O$ by 
$$\O_V=\Big\{x\in\O,\ \exists r_x>0\hbox{ such that }||V||_{L^{\frac n2,1}(\,B(x;r_x)\,)}<+\infty\Big\}.$$
One can show that $\O_V$ is an open set in $\O$.\\
Thus its complement $K_V=\O-\O_V$ is a compact included in $\O$.
\begin{defi}\ \\
The points $K_V$ are called the irregular points of $V$.
\end{defi}
\begin{rem}\ \\
The choice of $K_V$ can be modified  according to the application that one wants to do.
\end{rem}
If $V(x)=|x-a|^{-m},\ a\in\O$ and applying the first theorem, then
$$K_V==\begin{cases}\{a\}&if\ m\GEQ2,\\\emptyset&\hbox{otherwise}.\end{cases}$$
And as consequence of the above result, if $m>2$,  $A$ compact subset of $\O$
$$V(x)=\dist(x;A)^{-m} \hbox{ then } K_V=A.$$
\begin{coro}{\bf of Theorem \ref{t1} and Theorem \ref{t4}}\\
Under the same assumptions as for Theorem \ref{t4} and Theorem \ref{t1}, with $K=K_V$, then \\for $\Phi\in C^2_0(\OVO)$ the sequence $(\Phi_j)_j$ given in the proof of Theorem \ref{t4} say
$$\Phi_j=(1-\f_j)H\left(\dfrac2{\mu_j}\rho\right)\Phi$$
satisfies : 
 For all open set $\O_{V,0}$ relatively compact in $\O_V$ one has :
\begin{enumerate}
\item $\DST \Max_{\OVO_{V,0}}|\Phi_j(x)-\Phi(x)|\xj$
\item $|\nabla(\Phi_j-\Phi)|_{L^{n,1}(\O_{V,0})}\xj.$
\end{enumerate}
\end{coro}
\pr
Let $\O_{V,0}\subset\subset\O_V.$ Then according to Theorem \ref{t1} 
\begin{equation}\label{eq84}
\Max_{\OVO_{V,0}}|\f_j(x)|\xj
\end{equation}
\begin{equation}\label{eq85}
||\nabla\f_j||_{L^{n,1}(\O_{V,0})}\xj.
\end{equation}
On the other hand for $j\GEQ j_0$, we have $\O_{V,0}\subset\big\{\rho>\mu_j\big\}$. Therefore we have
\begin{equation}\label{eq86}
\Max_{\OVO_{V,0}}|\Phi_j(x)-\Phi(x)|\LEQ||\Phi||_\infty\Max_{\OVO_{V,0}}|\f_j(x)|
\end{equation}
\begin{equation}\label{eq87}
||\nabla(\Phi_j-\Phi)||_{L^{n,1}(\OVO_{V,0})}\LEQ||\nabla\Phi||_\infty\Max_{\OVO_{V,0}}|\f_j(x)|+||\Phi||_\infty\,||\nabla\f_j||_{L^{n,1}(\O_{V,0})}
\end{equation}

Relations (\ref{eq84}) to (\ref{eq87}) give the result.
\HF

As in \cite{DGRT}, we may add a transport term $U\cdot\nabla\f$ in the above equation (\ref{eq51}).

\begin{lem}\label{lA4}\ \\
Let $V$ a nonnegative potential  $K$ be a compact in $\O$ with $\Cap_\Vinf(K)=0$.\\
Consider $U\in L^{p,1}(\O)^n,\ p>n,\ w\in L^1(V)\cap L^q(\d^{-1}),\ q<p'$. Assume that $w\in L^{\frac n{n-2},\infty}(\O)$ if $n\GEQ3$ and $w\in L_{exp}(\O)$ if $n=2.$\\
Then, for all $\Phi\in C^2_0(\OVO)$, the sequence given in Theorem \ref{t4}, $\Phi_j=(1-\f_j)H\left(\dfrac2{\mu_j}\rho\right)\Phi$ satisfies
\begin{enumerate}
\item $\lim\DST\int_\O w\Delta\Phi_j\,dx=\int_\O w\Delta\Phi\,dx$
\item $\lim\DST\int_\O w\,\Phi_j \,Vdx=\int_\O\Phi\,w\,V\,dx$,
\item
$\DST\lim\int_\O wU\cdot\nabla\Phi_j\,dx=\int_\O wU\cdot\nabla\Phi\,dx$.
\end{enumerate}
\end{lem}
\pr
The two first statements are the consequence of the fact $w\in L^1(V)\cap L^1_0(\O)$, \ $(L^1(\d^{-1})\subset L^1_0(\O)$)  and the fact that $\Phi_j$ converges  weakly-strongly to $\Phi$ in the sense of the potential $V$. Moreover, we have a constant $c_{H\Phi}>0$
\begin{equation}\label{eq100}
\int_\O\Big|wU\cdot(\nabla\Phi_j-\nabla\Phi)\Big|dx\LEQ c_{H_\Phi}\left[\left\|\dfrac{\nabla\f_j}{\sqrt V}\right\|_\infty\int_\O|w|\,|U|\sqrt V dx+\int_\O|\f_j|\,|w|\,|U|dx\right].
\end{equation}
By H\"older, $w|U|$ and $|w|\,|U|\sqrt V$ are in $L^1(\O)$ since
$$\int_\O|w|\,|U| dx\LEQ\left\|\dfrac w\d\right\|_{L^q}||U||_{L^{q'}}<+\infty,\ q<n',\ \dfrac1q+\dfrac1{q'}=1$$
and
$$\int_\O|w|\,|U|\sqrt Vdx
\LEQ c||w||^{\frac12}_{L^1(V)}||w||^{\frac12}_{L^{\frac n{n-2},\infty}}
\cdot||U||_{L^{n,1}}
\ \si n\GEQ3.$$
Idem for $n=2$. Therefore, relation (\ref{eq100}) leads to statement 3. knowing 
$$\lim\f_j(x)=0,\ 0\LEQ\f_j\LEQ1\hbox{ and } \left\|\dfrac{\nabla\f_j}{\sqrt V}\right\|_\infty\xj.$$ 
\ \qquad\HF

\begin{theo}\label{t6}\ \\
Under the same assumption as for Lemma \ref{lA4}, if furthermore $w$ satisfies
\begin{equation}\label{eq101}
\int_\O w(-\Delta\Phi-U\cdot\nabla\Phi+ V\,\Phi)dx=0\quad\forall\,\Phi\in C^2_c(\O\BSL K)
\end{equation}
then, (\ref{eq101}) holds for all $\Phi\in C^2_0(\OVO)$,
and if $\p\O\in C^{1,1}$ and $\div(\ORA U)=0$ in $\calD'(\O)$ with $\ORA U\cdot \ORA\nu=0$ on $\p\O$ ($\nu$ exterior normal to $\p\O$) then
$$w\equiv0.$$
\end{theo}
\pr
If $\Phi\in C^2_0(\OVO)$ we have a sequence $\Phi_j$ in $C^2_c(\O\BSL K)$ such that, $\Phi_j$ satisfies the conclusion of Lemma \ref{lA4}. Thus, we have (\ref{eq101}) with $\Phi\in C^2_0(\OVO)$ as test function. To prove that $w\equiv 0$ we need to employ the following variant of Kato's inequality (see \cite{DGRT}).

\begin{theo}\label{tK}{\bf Variant of Kato's inequality}\ \\
Let $\OV u$ be in $W^{1,1}_{loc}(\O)\cap L^{n',\infty}(\O)$ with $\dfrac{\OV u}\d\in L^1(\O)$ and $\ORA U \in L^{n,1}(\O)^n$ with $\div(\ORA U)$ in $\calD'(\O),\ \ORA U\cdot\ORA\nu=0$ on $\p\O$.\\
Assume that $L\OV u=-\Delta\OV u+\div(\ORA U\OV u)\in L^1(\O;\d)$. Then for all $\phi\in C^2_0(\OVO), \ \phi\GEQ0$ one has
\begin{enumerate}
\item $\DST\int_\O \OV u_+L^*\phi\,dx\LEQ\int_\O\phi\sign_+(\OV u)L\OV u\,dx$
\item $\DST\int_\O|\OV u|L^*\phi\,dx\LEQ\int_\O\phi\sign(\u)L\OV u\,dx$,
\end{enumerate}
where $L^*\phi=-\Delta\phi-\ORA u\cdot\nabla\phi=-\Delta-\div(\ORA U\,\phi)$,

\end{theo}

According to equation (\ref{eq101}), $Lw=-Vw\in L^1(\O)$. Thus the above Kato's type inequality holds  and $$\forall\phi\in C^2_0(\OVO),\quad\phi\GEQ0:\int_\O|w|L^*\phi\LEQ-\int_\O\phi|w|Vdx\LEQ0.$$
Thus one has
$$\int_\O|w| L^*\phi=0\quad \forall\,\phi\in C^2_0(\OVO).$$

By density result the same equation holds
$$\forall\,\phi\in H^1_0(\O)\cap W^2L^{n,1}(\O),\ \phi\GEQ0.$$
\hbox{ Resolving }$L^*\phi=1$ we derive that $w\equiv0.$\HF

Next, we want to discuss some existence problem related to equation (\ref{eq101}). 
\\ 
We always assume that $U\in L^{p,1}(\O)^n,\ p>n$, $\div(U)=0$ in ${\mathcal D}'(\O)$, $U\cdot\nu=0$ on $\p\O$ and $\p\O\in C^{1,1}$.
\begin{theo}\label{t7}\ \\
Let $f$ be a bounded Radon measure in $\O$. Assume that $\Cap_\Vinf(K_V)=0$.\\
If $|f|(K_V)=0$ ($f$ does not charge the compact set $K_V$) then there exists an unique solution $u\in L^1(V)\cap L^1(\O;\d^{-1})$ such that
\begin{equation}\label{eq56}
\int_O u(-\Delta\f-U\cdot\nabla\f+V\f)dx=\int_\O\f df\qquad\forall\,\f\in C^2_0(\OVO).
\end{equation}
\end{theo}
\pr
The uniqueness is a consequence of Theorem \ref{t6}.\\
 For the existence, we first notice that the problem is linear, we may assume that $f\GEQ0$.  We shall set as usual
 $$M^1(\O)=\Big\{ f:\hbox{ bounded Radon measure on }\O\Big\},$$
 $$M^1(\O)=C_0(\OVO)',\quad C_0(\OVO)=\Big\{\f:\OVO\to\R\hbox{ continuous, $\f=0$ on }\p\O\Big\}$$
Let us introduce $V_j=\min(j;V)$ we have proved the following result in \cite{DGRT, DGR}. 

\begin{lem}\label{l7}\ \\
There exists $u_j\GEQ0,\ u_j\in W^1_0L^{n',\infty}(\O)$ such that
\begin{enumerate}

\item  $\forall\,\f\in H^1_0(\O)\cap W^2L^{n,1}(\O)$
\begin{equation}\label{eq57}
\int_\O u_j[-\Delta\f-U\cdot\nabla\f+V_j\f]dx=\int_\O\f df.
\end{equation}

\item There exists a constant $c_0$ independent of $j$ such that 
\begin{equation}\label{eq58}
||u_j||_{W^1_0L^{n',\infty}(\O)}+\int_\O V_ju_jdx\LEQ c_0||f||_{M^1(\O)}.
\end{equation}

\item In particular, there exist a function $u\GEQ 0$ and  a subsequence $u_j$ such that
\begin{enumerate}

\item $u_j\xji u(x)$ a.e in $\O$, strongly in $L^1(\O)$ and weakly in $W^1_0L^{n',\infty}(\O)$.

\item  
\begin{equation}\label{eq59}
||u||_{W^1_0L^{n',\infty}(\O)}+\int_\O Vu\,dx\LEQ c_0||f||_{M^1(\O)}.
\end{equation}
\end{enumerate}
\end{enumerate}
\end{lem}
{\bf Proof of Lemma \ref{l7}}\\
Since $f\in \Mu(\O)$, there is a sequence $f_k\in L^\infty_+(\O)$ such that 
$$||f_k||_{L^1(\O)}\LEQ||f||_{\Mu(\O)}\hbox{ and $f_k$ converges to $f$ weakly in $C_c(\O)'$}$$
(ie $\forall\,\f\in C_c(\O) \quad <f_k,\f>\xrightarrow[\   \ ]{}<f,\f>$.)

According to \cite{DGRT, DGR}, one has a function $u_{jk}\in W^1_0L^{n',\infty}(\O)$ satisfying, $\forall\,\f\in C_0^2(\OV\O)$
\begin{equation}\label{eq5000}
\int_\O u_{jk}\big[-\Delta\f-U\cdot\nabla\,\f+V_j\f\big]dx=\int_\O f_k\f dx
\end{equation}
and
\begin{equation}\label{eq6000}
||u_{jk}||_{W^1L^{n',\infty}(\O)}+\int_\O V_ju_{jk}dx\LEQ c_0||f_k||_{L^1(\O)}\LEQ c_0||f||_{\Mu(\O)}
\end{equation}
where $c_0$  is independent of $j$ and $k$ (in fact $c_0$ depends on $\O$ and $||U||_{L^{n,1}(\O)}$). More $u_{jk}\GEQ0$. Thus we have a subsequence  still denoted $(u_{jk})_k$ and a function $u_j\in W^1_0L^{n',\infty}(\O),\ u_j\GEQ0$ such
$\DST u_{jk}\mathop{\rightharpoonup}_{k\to+\infty} u_j\hbox{ weakly in $ W^1_0L^{n',\infty}(\O)$}$, strongly in $L^1(\O)$ and almost everywhere in $\O$.\\
Thus, we can pass easily to the limit in relations (\ref{eq5000}) and (\ref{eq6000}) to derive the part 1.) and 2.) of Lemma \ref{l7}. By the same reason as above, we have a subsequence still denoted $u_j$ and a function $u\GEQ0$ such that $u_j\rightharpoonup u$ weakly in $W^1_0L^{n',\infty}(\O)$ strongly in $L^1(\O)$, almost everywhere in $\O$. From relation (\ref{eq6000}) using among other Fatou's lemma, we have relation (\ref{eq59}).\HF
\begin{lem}\label{l8}\ \\
Let $\f\in W^1_0L^{n,1}(\O)$ with support$(\f)\cap K_V=\emptyset$.\\
Then
$$\lim_{j\to+\infty}\int_\O|u_jV_j-uV|\,|\f|dx=0.$$
\end{lem}
\pr
Let $\f$ be in $W^1_0L^{n,1}(\O)$ with support$(\f)\cap K_V=\emptyset$.\\

Thus 
$V\f\in L^{\frac n2,1}(\O)$ and support$(V\f)\subset\subset \O\BSL K_V=\O_V$, \\
 since support$(\f)\subset\Big\{x:\dist(x;K)>\eta\Big\}$ for some $\eta>0$. We have :
\begin{eqnarray}\label{eq60}
||u_j-u||_{L_{exp}(\O)}&\LEQ& c_N||u_j-u||_{W^1_0L^{2,\infty}(\O)}\ \si n=2\NN\\
||u_j-u||_{L^{\frac n{n-2},\infty}(\O)}&\LEQ& c_N||u_j-u||_{W^1_0L^{n',\infty}(\O)}\ \si n\GEQ3.
\end{eqnarray}
Therefore, applying H\"older's inequality, we have a constant $c_6>0$ (independent of $u_j,\ u,\ V_j$) such that for any measurable subset $E\subset\O$
\begin{equation}\label{eq61}
\int_E   V\,|\f|\,|u_j-u|dx\LEQ c_6||V\f||_{L^{\frac n2,1}(E)}\xrightarrow[|E|\to0]{}0.
\end{equation}
Therefore, using the Egoroff's theorem or Vitali's theorem one has
\begin{equation}\label{eq62}
\lim_{j\to\infty}\int_\O V\,|\f|\,|u_j-u|dx=0.
\end{equation}
Since we have
\begin{equation}\label{eq63}
\lim_{j\to\infty}\int_\O |u|\,|V_j-V|\,|\f| dx=0.
\end{equation}
Finally we have
\begin{equation}\label{eq64}
\lim_{j\to\infty}\int_\O|u_jV_j-uV|\,|\f|dx=0.
\end{equation}
\ \qquad\HF

\begin{lem}\label{l9}\ \\
The function $u$ found in the preceding Lemma \ref{l7}, satisfies, $\forall\,\f\in C^2_c(\O_V)$
$$\int_\O u(-\Delta\f-U\cdot\nabla\,\f+V\f)dx=\int_\O\f df.$$
\end{lem}

\pr
Let $\f$ be in $C^2_c(\O_V)$ then,
\begin{equation}\label{eq65}
\int_\O u_j[-\Delta\f-U\cdot\nabla\,\f+V_j\f]dx=\int_\O\f\,df
\end{equation}
and
$$\lim_{j\to+\infty}\int_\O|u_jV_j-uV|\,|\f|dx=0\hbox{ (since support}(\f)\cap K_V=\emptyset).$$
Thus we may pass to the limit in relation (\ref{eq65}).
\HF

\begin{lem}\label{l10}\ \\
If $\Cap_\Vinf(K_V)=0$ and $|f|(K_V)=0$ then, 
$$\hbox{the function $u$ given in Lemma \ref{l9} satisfies  relation (\ref{eq56}) }$$and $$u \in L^1(V)\cap W^1_0L^{n',\infty}(\O)\subset L^1(V)\cap L^q\left(\O;\dfrac1\d\right)\quad q<n'.$$
\end{lem}
\pr
Let $\Phi$ be in $C^2_0(\OVO)$. From our assumption we have a sequence $\Phi_j\in C^2_c(\O_V)$ such that :
\begin{enumerate}

\item $\Phi_j$ converges weakly to $\Phi$ in the sense of potential $V$,

\item $\Phi_j(x)\xji\Phi(x)$ for all $x\in\O-K_V=\O_V$ (see Corollary 7.2 of Theorem \ref{t4} and Theorem \ref{t1}), $\Phi_j=0$ in the neighborhood of $K_V$
\end{enumerate}
Since
\begin{equation}\label{eq70}
\int_\O u(-\Delta\Phi_j-U\cdot\nabla\,\f_j+V\Phi_j)dx=\int_\O\Phi_j df,
\end{equation}
We pass to the limit since $u\in L^1(V)\cap  W^1_0L^{n',\infty}\O\subset L^1(V)\cap L^1_0(\O)$ in the first integral and in the second integral using Lebesgue dominated theorem to derive
\begin{eqnarray}
\lim_{j\to\infty}\int_\O u(-\Delta\Phi_j+V\Phi_j)dx&=&\int_\O u(-\Delta\Phi+V\Phi)dx\label{eq71}\\
\lim_{j\to\infty}\int_\O\Phi_jdf&=&\int_{\O\BSL K_V}\Phi df.\label{eq72}
\end{eqnarray}
Applying  Lemma \ref{l7}, we have 
\begin{equation}\label{eq167}
\lim_{j\to\infty}\int_\O u\,U\cdot\nabla\Phi_jdx=\int_\O u\,U\cdot\nabla\Phi dx.
\end{equation}
But $|f|(K_V)=0$ so we have 
\begin{equation}\label{eq73}
\int_{\O\BSL K_V}\Phi df=\int_\O\Phi df.
\end{equation}
From relations (\ref{eq70}) to(\ref{eq73}), we derive 
\begin{equation}\label{eq74}
\int_\O u(-\Delta\Phi-U\cdot\nabla\,\Phi+V\Phi)dx=\int_\O\Phi df.
\end{equation}
\ \qquad\HF

\newpage

For the converse, we will first prove

\begin{theo}\label{t8}\ \\
Assume that $\Cap_\Vinf(K_V)=0,\ f=\mu_a$ , the Dirac measure at $a\in\O$.\\
If $a\in K_V$ then there is no solution  $u\in L^1(V)\cap L^1_0(\O)$ of
\begin{equation}\label{eq75}
\int_\O u[-\Delta\f-U\cdot\nabla\,\f+V\f]dx=\f(a)\qquad\forall\f\in C^2_0(\OVO).
\end{equation}
\end{theo}

\pr
If  there was a solution, then, $\forall \f\in C^2_c(\O\BSL K_V)$ we have
$$\int_\O u[-\Delta\f-U\cdot\nabla\,\f+V\f]dx=0.$$
But $\Cap_\Vinf(K_V)=0$ thus the same equation holds for all $\f\in C_0^2(\OVO)$ which implies that $\forall\,\f\in C_0^2(\OVO),\ \f(a)=0$. This is impossible.\HF

One can generalize Theorem \ref{t8} as follow

\begin{theo}\label{t9}\ \\
Assume that $\Cap_\Vinf(K_V)=0$. Let $f$ be a bounded Radon measure such that\\ $G=$support$(f)\cap K_V$ is an isolate subset of support$(f)$, ie. there exists an open set $\omega$ such that $\omega\cap $support$(f)=G$.\\
 
\begin{equation}\label{eq76}
\si |f|(K_V)>0\hbox{ then there is no solution of (\ref{eq74}).}
\end{equation}
\end{theo}
\pr
If $|f|(K_V)>0$, then $G$ is an isolate subset of support of $f$, therefore, we can consider $\theta\in C^\infty_c(\O)$ such that $\theta=1$ on $G$, support$\theta\subset \omega$.\\ We write $f=\theta f+(1-\theta) f$ so that measure $f_1=(1-\theta)f$ does not charge $K_V$. By the preceding result, we have $u_1\in L^1(V)\cap L^1_0(\O)$ such that
$$\int_\O u_1[-\Delta\f-U\cdot\nabla\,\f+V\f]dx=\int_\O\f df_1\quad \forall\f\in C^2_0(\OVO).$$
Assume that we have a solution $u$ of (\ref{eq74}) so, $w=u-u_1$ is a solution of 
$$\int_\O w[-\Delta\f-U\cdot\nabla\,\f+V\f]dx=<\theta f,\f>,\qquad\forall\,\f\in C^2_0(\OVO).$$
In particular
$$\int_\O w[-\Delta\f-U\cdot\nabla\,\f+V\f]=0\qquad\forall\,\f\in C^2_c(\O\BSL K_V).$$
Applying Theorem \ref{t6}, we deduce that $w=0$ say $u=u_1$ which mean $f=(1-\theta )f:\theta f\equiv0$ this is a contradiction with the fact that $|f|(K_V)>0$.
\HF
\begin{theo}\label{t10}\ \\
Assume that $\Cap_\Vinf(K_V)=0$, $f\in M^1(\O)$ such that $G={\rm support\,}(f)\cap K$ is an isolate subset of ${\rm support\,}(f)$. \\
Then one has a solution $u\in L^1(V)\cap L^1_0(\O)$ of 
$$\int_\O u[-\Delta\f-U\cdot\nabla\,\f+V\f]=\int_\O\f df\ \ \forall\f\in C^2_0(\OVO)\hbox{ if and only if } |f|(K_V)=0.$$
\end{theo}
$$\underline{\ \qquad\ \ \qquad}$$
\ \\
After submitting this work, we have received the paper \cite{OP} where a similar result as for this last theorem is given but only for solution in $W^{1,1}_0(\O)\cap L^1(V)$ which is strictly included in $L^1_0(\O)\cap L^1(V)$.\\
More, our proofs are totally different.

\ \\

{\bf Acknowledgment}\\
This work was initiated partly in December 18$^{th}$, 2017 when the author  attended to the conference {"\it Nonlinear Partial Equation and Mathematical Analysis"}\\
He would like to thank all the participants  and the organizers for their warm hospitality and invitation. A special thanks to Prof. D\'iaz Ildefonso whose friendship is a constant encouragement and an inspiration for him.

He wants also to thank the anonymous referees for  reading carefully this manuscript.

\end{document}